\documentclass[preprint,12pt]{elsarticle}

\usepackage{amssymb}

\newtheorem{thm}{Theorem}[section]
\newtheorem{cor}[thm]{Corollary}
\newtheorem{prp}[thm]{Proposition}
\newtheorem{lem}[thm]{Lemma}
\newproof{pf}{Proof}
\newdefinition{df}[thm]{Definition}
\newdefinition{exm}[thm]{Example}

\newcommand{\CC}{{\mathbb{C}}}

\newcommand{\NN}{{\mathbb{N}}}

\newcommand{\RR}{{\mathbb{R}}}

\newcommand{\cG}{{\cal G}}
\newcommand{\cS}{{\cal S}}
\newcommand{\cC}{{\cal C}}
\newcommand{\ext}{{\rm Ext \,}}

\begin{document}

\begin{frontmatter}

\title{On semigroups of matrices with nonnegative diagonals\tnoteref{label1}}

\tnotetext[label1]{The paper will appear in Linear Algebra and its Applications. \\
The authors were supported by the Slovenian Research Agency. \\ }

\author{Grega Cigler}
\ead{gregor.cigler@fmf.uni-lj.si}

\author{Roman Drnov\v{s}ek\corref{cor1}}
\ead{roman.drnovsek@fmf.uni-lj.si}

\cortext[cor1]{Corresponding author}

\address{Department of Mathematics, Faculty of Mathematics and Physics, University of Ljubljana, 
Jadranska 19, SI-1000 Ljubljana, Slovenia}

\baselineskip 6mm

\begin{abstract}
We give a short proof of a recent result by Bernik, Mastnak, and Radjavi,
stating that an irreducible group of complex matrices with nonnegative
diagonal entries is diagonally similar to a group of  nonnegative monomial matrices. 
We also explore the problem when an irreducible matrix semigroup in which each member is diagonally similar to
a nonnegative matrix is diagonally similar to a semigroup of nonnegative matrices. 
\end{abstract}


\end{frontmatter}

\baselineskip 6.8mm



\section{Introduction}

Multiplicative semigroups of matrices with nonnegative diagonal entries have been studied in the papers 
\cite{BMR} and \cite{SWG}. Their authors considered the general question under which additional assumptions
such a semigroup is simultaneously similar to a semigroup of nonnegative matrices.
The main result of \cite{BMR} is that every irreducible group of complex matrices with nonnegative
diagonal entries is diagonally similar to a group of  nonnegative monomial matrices.
In Section 2 we give a short proof of this result. Our proof is more geometric and less group-theoretic than the proof in \cite{BMR}.
Multiple authors of the paper \cite{SWG} provided several examples showing that it is impossible to extend this result 
from groups to semigroups. So, to obtain similarity to a semigroup of nonnegative matrices, stronger assumptions on a given semigroup 
must be imposed. In Section 3 we explore the problem when an irreducible matrix semigroup in which each member is diagonally similar to
a nonnegative matrix is necessarily diagonally similar to a semigroup of nonnegative matrices. 

We now recall some definitions and basic facts. 
The set of all nonnegative real numbers is denoted by $\RR_+$. 
A convex set $K \subseteq \RR^n$ is said to be a {\it cone} if $r K \subseteq K$ for all $r \in \RR_+$.  
A cone $K  \subseteq \RR^n$ is {\it proper} if it is closed, {\it pointed} ($K \cap (-K) = \{0\}$), and {\it solid} (the interior of $K$ is nonempty). 
The most natural example of a proper cone is the {\it nonnegative orthant} $\RR_+^n$.
A cone $K \subseteq \RR^n$ is {\it reproducing} if $K - K = \RR^n$. It is well-known that a closed cone is solid if and only if it is reproducing.

Let $K$ be a closed cone in $\RR^n$. A vector $x \in K$ is an {\it extremal vector}
of $K$ if $y \in K$ and $x-y \in K$ imply that $y$ is a nonnegative multiple of $x$.  
By $\ext (K)$ we denote the set of all extremal vectors of $K$.
By the Krein-Milman theorem, $K$ is the convex hull of $\ext (K)$.
The angle $\phi \in [0, \pi]$ between non-zero vectors $x$, $y \in \RR^n$ is
determined by the equality \ $x^T y = \|x\| \, \|y\| \, \cos \phi$.

If $F$ is a subset of complex numbers, then $M_n(F)$ denotes the set of all $n \times n$ matrices with entries in $F$.
If $\cC \subseteq M_n(\CC)$ is a collection of complex matrices, then $\overline{\cC}$ denotes its closure in the Euclidean topology, 
and $\RR_+ \cC$ denotes its {\it homogenization}, i.e., $\RR_+ \cC = \{r C: r \in \RR_+, C \in \cC\}$.
We say that a matrix has a {\it nonnegative diagonal} if all of its diagonal entries are nonnegative.
A matrix is called {\it monomial} if it has the same nonzero pattern as a permutation matrix, i.e., 
there is exactly one nonzero entry in each row and in each column.

A collection $\cC \subseteq M_n(\CC)$ (where $n \ge 2$) is {\it reducible} if there exists a common invariant subspace other than 
the trivial ones $\{0\}$ and $\CC^n$, or equivalently, there exists an invertible matrix $S \in M_n(\CC)$ such that the collection 
$S \cC S^{-1}$ has a block upper-triangular form; otherwise, the collection $\cC$ is said to be {\it irreducible}.
If the matrix $S$ can be chosen to be a permutation matrix, then the collection $\cC$ is said to be {\it decomposable};
otherwise, it is called {\it indecomposable} (or {\it ideal-irreducible}).

\section{Groups of matrices with nonnegative diagonals}

The study of semigroups of matrices having nonnegative diagonals was initiated by the authors of \cite{BMR}.
They started their discussion by the following result (see \cite[Theorem 4.1]{BMR}).

\begin{thm}
\label{rank-one}
Let $\cS \subseteq M_n(\CC)$ be an irreducible semigroup of matrices of rank at most one having nonnegative diagonals.
If  $\overline{\RR_+ \cS} = \cS$, then, after a diagonal similarity, $\cS = X Y^T$ for some subsets $X$ and $Y$ of $\RR_+^n$ 
each of which spans $\CC^n$. 
\end{thm}

Using the Haar measure one can prove the following assertion (see \cite[Proposition 4.3]{BMR}).

\begin{prp}
\label{positive-valued}
Let $\cS \subseteq M_n(\CC)$ be an irreducible semigroup of matrices. 
Suppose that $\overline{\RR_+ \cS} = \cS$ and that there exists a non-zero functional $\varphi: M_n(\CC) \to \CC$
such that $\varphi(S) \in \RR_+$ for all $S \in \cS$. Then $\cS$ has members of rank one. 
\end{prp}

The following theorem is the main result of \cite[Theorem 5.5]{BMR}. 
We provide a short proof that is more geometric and less group-theoretic than the original one.

\begin{thm}
\label{positivediaggroup}

If  $\cG \subset M_n(\CC)$ is an irreducible group of matrices with nonnegative diagonals, then, up
to a diagonal similarity, $\cG$ is a group in $M_n(\RR_+)$. 
Therefore, each member of the group $\cG$ is a nonnegative monomial matrix.
\end{thm}

\begin{pf} 
With no loss of generality we may assume that $t G \in \cG$ for all $t > 0$ and $G
\in \cG$. Let $\cS =  \overline{\cG}$. Applying Proposition \ref{positive-valued} for the trace
functional, we conclude that $\cS$ contains elements of rank one. The semigroup ideal $\cS_1$ of all elements of
rank at most one in $\cS$ is irreducible (see \cite{RR}). By Theorem \ref{rank-one}, we can assume that,
after a diagonal similarity, $\cS_1 = X Y^T$ for some subsets $X$ and $Y$ of $\RR_+^n$ each of which spans
$\CC^n$. We can also assume that $\RR_+ X = X$ and $\RR_+ Y = Y$. 
The cone $\widehat{X}$ generated by $X$ is closed, and it is invariant under any $S
\in \cS$, since $ (Sx) y^T = S (x y^T) \in \cS_1$ for every $x \in X$ and $y \in Y$. 
Similarly, it follows from $x (S^T y)^T = (x y^T) S \in \cS_1$ that $Y$ is invariant under $S^T$.
The dual cone  
$$ Y^d = \{ z \in \RR^n :  z^T y \ge 0 \textrm{  for all  } y \in Y \} $$
of the set $Y$ obviously contains $\RR_{+}^n$, and it is invariant under any $S \in \cS$, as 
$(Sz)^T y = z (S^T y) \ge 0$ for all $y \in Y$ and $z \in Y^d$.  
It follows that every $G \in \cG$ is a bijective mapping on both $\widehat{X}$ and
$Y^d$, implying that every $G \in \cG$ maps $\ext(\widehat{X})$ to itself, and the same
holds for the cone $Y^d$.
We want to show that the inclusions $\widehat{X} \subseteq \RR_+^n \subseteq  Y^d$
are in fact equalities.

Assume, if possible, that $\widehat{X} \neq Y^d$. Then there exists a unit vector 
$x \in X \setminus Y^d$ which is extremal for the cone $\widehat{X}$.
Since the cone $Y^d$ is closed, the distance between $x$ and $Y^d$ is strictly positive.
It follows that there is a number $\phi \in (0, \pi/2)$
such that, for each $z \in \ext(Y^d)$, the angle between $z$ and $x$ is at least $\phi$. 
Since $x \in X$ and the set $Y$ is spanning, there is a vector $y \in Y$ such that
$P = x y^T \in \cS$ with $y^T x > 0$. 
We can assume that $y^T x = 1$, so that $P x = x$.
Choose any $\epsilon > 0$. 
Since  $\cS =  \overline{\cG}$, there is a matrix $G \in \cG$ such that $\|G - P\| < \epsilon$. 
Now, for any $z \in \ext(Y^d)$ with norm $1$, we have 
$$ \epsilon^2 > \|G z - P z\|^2 = \| G z - (y^T z) x \|^2 = 
   \|G z\|^2 + (y^T z)^2 - 2 (y^T z) \|G z\| \cos \phi_z , $$
where $\phi_z$ is the angle between the vector $x$ and the vector $G z \in \ext(Y^d)$.
Since $y^T z \in \RR_+$ and $\phi_z \geq \phi$, we conclude that 
$$ \epsilon^2 > \|G z\|^2 + (y^T z)^2 - 2 (y^T z) \|Gz\| \cos \phi = 
   (y^T z - \|G z\| \cos \phi)^2 + \|G z\|^2 \sin^2 \phi . $$
It follows that 
$$ \|G z\| \sin \phi < \epsilon \ \  \textrm{and} \ \ \ \left|y^T z - \|G z\| \cos \phi \right| < \epsilon , $$
and so 
$$ 0 \le y^T z < \epsilon + \|G z\| \cos \phi < \epsilon +  \frac{\epsilon}{\sin \phi} \cos \phi . $$
Since $\epsilon>0$ is arbitrary, we obtain that $y^T z = 0$ for all vectors $z \in \ext(Y^d)$,
implying that $y = 0$. 
This contradiction completes the proof of the equality $\widehat{X} = Y^d = \RR_+^n$.
Consequently,  the inclusion $\cG \subset M_n(\RR_+)$ holds, as asserted.

Since the map associated to any matrix $G \in \cG$ maps $\ext(\RR_+^n)$ to itself 
and it is invertible, the matrix $G$ must be monomial, and so the proof is complete.
\qed \end{pf}

\section{Semigroups of matrices diagonally similar to nonnegative ones}

Let $\cS \subseteq M_n(\CC)$ be a semigroup in which each member $A\in \cS$ is diagonally similar to
a nonnegative matrix. In this section we are looking for additional assumptions under which 
the whole semigroup $\cS$ is diagonally similar to a semigroup of nonnegative matrices. 
We first show that it does not suffice to assume that the semigroup $\cS$ is indecomposable.

\begin{exm} 
Define $n \times n$ matrices $A=a a^T$ and $B=b b^T$, where $n \ge 2$, $a=[1,1,\ldots,1]^T$ and $b=[1,1,\ldots,1,1-n]^T$. 
Then every nonzero member of  the semigroup $\cS$ generated by $A$ and $B$ is an indecomposable matrix of rank one that is diagonally similar to a nonnegative matrix. 
However, the whole semigroup $\cS$ is not diagonally similar to a semigroup of nonnegative matrices.
\end{exm}

\begin{pf}
Note that $A^k=n^{k-1} A$ and $B^k=(n(n-1))^{k-1} B$ for all $k\in\NN$, while $AB=BA=0$. Therefore, $\cS$ is contained 
in the semigroup $\RR_+ A \cup \RR_+ B$.  
If $D$ is the diagonal matrix with diagonal $(1,1,\ldots,1,-1)$, then the matrix $DBD^{-1}$ is nonnegative, and therefore
each matrix from $\cS$ is diagonally similar to a nonnegative matrix. Since the matrices $A$ and $B$ are indecomposable, 
every nonzero member of  $\cS$ is indecomposable as well. It is easy to verify that the whole semigroup $\cS$ is not 
diagonally similar to a semigroup of nonnegative matrices.
\qed
\end{pf}

In the rest of the paper we explore the case when the semigroup $\cS$ is irreducible.
We first show that, with no loss of generality, we may assume that $\cS$ is a closed set.

\begin{lem}\label{closurelemma}
Let $\cC\subset M_n(\CC)$ be a collection in which each member $A\in \cC$ is diagonally similar to a nonnegative matrix.
Then the closure $\overline{\RR_+\cC}$ also consists of matrices which are diagonally similar to 
nonnegative matrices.
\end{lem}

\begin{pf}
Clearly, we may assume that $\RR_+ \cC = \cC$. If $A \in \overline{\cC}$, then there is a sequence $\{A_k\}_{k \in \NN}$ in $\cC$ 
converging to the matrix $A$. For each $k\in\NN$, let $D_k$ be a diagonal matrix such that $D_k A_k D_k^{-1}$ is a nonnegative matrix. 
We may assume that each diagonal entry of $D_k$ has absolute value one. 
Since the sequence $\{D_k\}_{k \in \NN}$ is bounded, it has a convergent subsequence $\{D_{k_m}\}_{m \in \NN}$
converging to some diagonal matrix $D$. Since  $DAD^{-1} =\lim_{m\to\infty} D_{k_m}A_{k_m}D_{k_m}^{-1}$,
the matrix $DAD^{-1}$ is nonnegative, and so $A$ is also diagonally similar to a nonnegative matrix. This completes the proof.
\qed
\end{pf}

We continue with a reduction of the problem to the real setting.

\begin{lem}\label{backtoreality} 
Let $\cS = \overline{\RR_+\cS} \subseteq M_n(\CC)$ be an irreducible semigroup such that each member $A\in \cS$ is diagonally similar to a nonnegative matrix.
Then there exists an invertible diagonal matrix $D \in M_n(\CC)$ such that the semigroup $D \cS D^{-1}$ consists of real matrices, 
and there exist two sets $X,Y \subseteq \RR_+^n$, each of which spans $\CC^n$, such that 
$$ D\cS_1D^{-1}=(D\cS D^{-1})_1=XY^T ,$$
where $\cS_1$ is the ideal of $\cS$ consisting of members of rank at most one.
Furthermore, the subcone of $\RR_+^{n}$ generated by $X$ is a proper cone invariant under 
every member of $\cS$.
\end{lem}

\begin{pf}
Our assumption implies in particular that all diagonal elements of any member of $\cS$ must be nonnegative.
By Proposition \ref{positive-valued}, the ideal $\cS_1$ of all members of $\cS$ with rank at most one is nonzero.
Since $\cS$ is an irreducible semigroup, it is also necessarily irreducible (see \cite{RR}).
Then by Theorem \ref{rank-one} we can find an invertible diagonal matrix $D$ and two sets $X,Y\subset \RR_+^n$, 
each of which spans $\CC^n$, such that $D \cS_1 D^{-1} = X Y^T$. 
As we are interested in diagonal similarities, we can assume that $D$ is the identity, so that $\cS _1=XY^T$.
To prove the inclusion $\cS \subset M_n(\RR)$, pick any $A\in \cS$ and $x\in X$. Since for any nonzero vector $y\in Y$ 
the matrix $A(xy^T)=(Ax)y^T$ belongs to $\cS_1$, we conclude that $Ax \in X \subseteq \RR_+^n$. It follows that the cone 
of $\RR_+^{n}$ generated by $X$ is a proper cone invariant under $A$.
Since the set $X$ spans $\CC^n$, it follows that $A (\RR^n) \subseteq \RR^n$, and therefore
$A \in M_n(\RR)$. This completes the proof.
\qed
\end{pf}

From now on we consider real matrices. If a real matrix $A$ is diagonally similar to a nonnegative matrix via diagonal matrix $D$, 
we clearly may assume that each diagonal entry of $D$ is either $1$ or $-1$. In this case we say that $D$ is a {\it $\pm 1$-diagonal} matrix.

\begin{lem}
\label{indecompobstruction}
Let $A \in M_n(\RR)$ be an indecomposable matrix and $D$ a $\pm 1$-diagonal matrix such that $A'=DAD$ 
is a nonnegative matrix. If there exists a proper cone $K$ such that $A(K)\subseteq K$ and $K\subseteq \RR^n_+$, then 
$D=\pm I$ and $A$ itself is a nonnegative matrix.
\end{lem}

\begin{pf}
By the Perron-Frobenius Theorem, the spectral radius $\rho(A')=\rho(A)$ of the indecomposable matrix $A'$
is a simple eigenvalue having exactly one (up to a scalar multiplication) strictly positive eigenvector $e$.
On the other hand, since the proper cone $K$ is invariant under $A$, the extension of the Perron-Frobenius Theorem 
(see \cite[Theorem 3.2]{BP}) ensures that there is a non-zero vector  $x \in K$ such that $A x = \rho(A) x$. 
However, $A' D x= D A x =\rho(A) D x$, and so the vectors $Dx$ and $e$ are collinear.
It follows that either $De$ or $-De$ belongs to $K \subseteq \RR_+^{n}$, 
and this implies that $D=\pm I$ and $A$ itself is a nonnegative matrix.
\qed
\end{pf}

The following simple example shows that in Lemma \ref{indecompobstruction} we cannot omit the assumption 
that the cone $K$ is proper.

\begin{exm}
Let $n \ge 2$, $a=[1,1, \ldots,1,1-n]^T$ and  $K=\RR_+ [1,1, \ldots,1]^T$. The matrix $A=a a^T$ is indecomposable, and the cone $K$ 
is invariant under $A$, while $DAD$ is a nonnegative matrix for the diagonal matrix $D$ with diagonal $(1,1, \ldots,1,-1)$.
\qed
\end{exm}

For $n\ge 2$ we say that a matrix $A\in M_n(\RR)$ is $1$-{\it decomposable} if there is a permutation matrix $P$ such that
$$ PAP^T=\left[\matrix {A_1  &   B    \cr
                          0    &   A_2     
}\right], $$
where each of $A_1$ and $A_2$ is either an indecomposable (square) matrix or a $1 \times 1$ block.

The following assertion is crucial for the proof of the main result.

\begin{prp}
\label{obstruction}
Let $A \in M_n(\RR)$ be a $1$-decomposable matrix that is diagonally similar to a nonnegative matrix. Let $K$ and $L$ be 
proper cones of $\RR^n_+$ that are invariant under $A$ and $A^T$, respectively. Then $A$ is a nonnegative matrix.
\end{prp}

\begin{pf}
Let $P$ be a permutation matrix such that the matrix $PAP^T$ has the block form
$$ PAP^T=\left[\matrix {A_1  &   B    \cr
                          0    &   A_2     
}\right] $$
with respect to the decomposition $\RR^n=\RR^k\oplus\RR^l$, where $1 \le k < n$, $l = n - k$, and 
each of $A_1$ and $A_2$ is either an indecomposable (square) matrix or a $1 \times 1$ block.
We first prove that the diagonal blocks $A_1$ and $A_2$ are nonnegative matrices.
If $DAD$ is a nonnegative matrix for a suitable $\pm 1$-diagonal matrix $D$, then $E=PDP^T$ is a $\pm 1$-diagonal matrix
such that $E(PAP^T) E$ is a nonnegative matrix. It follows that matrix $PAP^T$ satisfies our assumptions provided that the cones $K$ and $L$ are replaced by the cones $P(K)$ and $P(L)$. We can therefore assume that $A$ itself is of the block form
$$A=\left[\matrix {A_1  &   B    \cr
                   0    &   A_2     
}\right] \ . $$
Let $\Pi_1:\RR^n\to\RR^k$ and $\Pi_2:\RR^n\to\RR^l$ be the corresponding projections, and let $C\subseteq \RR^n_+$ be a proper cone.
As $C\subseteq\Pi_1(C)+\Pi_2(C)$ and $\Pi_1(C)$ contains at most $k$ linearly independent vectors, it follows that $\Pi_2(C)$ contains
at least $n-k=l$ linearly independent vectors. Consequently, $\Pi_2(C)$ contains exactly $l$ linearly independent vectors, so that $\Pi_2(C)$
is a generating cone of $\RR^l$.  Similarly, $\Pi_1(C)$ is a generating cone of $\RR^k$. 
Since $C\subseteq \RR^n_+$, both $\Pi_1(C)$ and $\Pi_2(C)$ are pointed and therefore proper cones. 
Assume now that the cone $C$ is invariant under $A$.  If  $x_2\in \Pi_2(C)$, then $x_2=\Pi_2(x)$ for some $x\in C$, and so $A_2(x_2)=A_2(\Pi_2(x))=\Pi_2(Ax)\in\Pi_2(C)$,
since $A(C)\subseteq C$. Therefore, the cone $\Pi_2(C)$ is invariant under $A_2$.
This means that $\Pi_2(K)\subseteq \RR^l_+$ is a proper cone invariant under $A_2$. Since the indecomposable matrix $A_2$ is diagonally similar to a nonnegative matrix, we can apply Lemma \ref{indecompobstruction} to conclude that $A_2$ is a nonnegative matrix.

In order to show that  $A_1$ is also a nonnegative matrix, we consider the transposed matrix $A^T$. The proper cone $L\subseteq \RR^n_+$ is invariant under $A^T$. 
Then the cone $\Pi_1(L)$ is a proper cone invariant under $A_1^T$. Since $A_1$ is indecomposable, $A_1^T$ is indecomposable and again by Lemma 
\ref{indecompobstruction} we conclude that $A_1$ must be a nonnegative matrix.

It remains to prove that the block $B$ is nonnegative.
Suppose to the contrary that $B$ has some strictly negative entries. 
If $D=D_1\oplus D_2$ is a $\pm 1$-diagonal matrix such that $DAD$ is a nonnegative matrix, 
then $D_iA_iD_i$ for $i=1,2$ and $D_1BD_2$ are nonnegative matrices. Using Lemma \ref{indecompobstruction} we conclude
that $D_i=\pm I$ for $i=1,2$ and $D_1BD_2=\pm B$. Since $B$ contains some strictly negative entries,  the matrix $-B$ must be nonnegative.
Since we can add the identity matrix to the matrix $A$, without loss of generality we can assume that the matrices $A_1$ and $A_2$ 
are both primitive, i.e., the spectral radius $\rho(A_i)$ is the only point in the peripheral spectrum of $A_i$, $i=1,2$. 
For $k\in\NN$ we have
$$A^k=\left[\matrix {A_1^k  &   B_k\cr
                   0    &   A_2^k     
}\right] ,$$
where 
$$B_k=\sum_{l=0}^{k-1}A_1^{k-1-l}BA_2^{l}.$$
If we multiply the matrix $A$ by a suitable positive scalar, we can assume that 
$\rho(A) = \max\{\rho(A_1), \rho(A_2)\}=1$.  We must consider the following three cases:

(1) $\rho(A_1)=\rho(A_2)=1$: By Perron-Frobenius theory, the limits
$$ \lim_{k\to\infty} A_1^k=E_1 \textrm{      and      }  \lim_{k\to\infty} A_2^k=E_2 $$ 
are strictly positive idempotents of rank $1$. 
In particular, there is a constant $C>0$ such that $\|A_1^k\|,\|A_2^k\|\le C$ for all $k\in \NN$. 
Then we have, for any $m\in\NN$,   
$$ \|B_{4m}\| = \left\| \sum_{l=0}^{4m-1} A_1^{4m-1-l}BA_2^{l} \right\| \le 
   \sum_{l=0}^{4m-1}\|A_1^{4m-1-l}\|\|B\|\|A_2^{l}\|\le 4m \, C^2 \|B\| , $$
and so the sequence $\{\frac{1}{4m} B_{4m}\}_{m \in \NN}$ is bounded. It follows that some subsequence $\{\frac 1{4m_k}A^{4m_k}\}_{k \in \NN}$ 
of the sequence $\{\frac{1}{4m} A_{4m}\}_{m \in \NN}$ converges to the matrix of the form
$$ A_\infty= \lim_{k\to \infty}\frac 1{4m_k}A^{4m_k}=\left[\matrix {0 &   B_\infty\cr
                                                    0    &   0     
}\right].$$
Choose $m \in \NN$ such that $\frac 12 E_i \le A_i^{l}$ for $i=1,2$ and all $l \ge m$. As $-B$ is a nonnegative matrix, we obtain that
$A_1^{4m-1-l}BA_2^{l} \le \frac 14 E_1BE_2$ for all $l = m, m+1, m+2, \ldots, 3m-1$. 
Since the matrices $- A_1^{4m-1-l} B A_2^{l}$ are nonnegative, we have 
$$ B_{4m} = \sum_{l=0}^{4m-1}A_1^{4m-1-l}BA_2^{l} \le \sum_{l=m}^{3m-1}A_1^{4m-1-l}BA_2^{l}\le \frac 14\sum_{l=m}^{3m-1}E_1BE_2 . $$
It follows that
$$ B_\infty \le\lim_{m\to\infty} \frac 1{4m}\left(\frac 14\sum_{l=m}^{3m-1}E_1BE_2\right)=\frac 18 E_1BE_2 , $$
and so $B_\infty$ is a matrix with some strictly negative entries.
Therefore, there is a strictly positive vector $e \in K$ such that the vector $A_\infty e$ is not in $\RR^n_+$.
As the cone $K$ is closed and invariant under all powers of $A$, it has to be invariant under $A_\infty$, so that 
$A_\infty e \in K \subseteq \RR^n_+$. This contradiction completes the proof in this case.

(2) $1=\rho(A_1)>\rho(A_2)$: As before, the limit $\lim_{k\to\infty} A_1^k=E_1$ is a strictly positive idempotent of rank $1$. 
Since $L \subseteq \RR^n_+$ is a proper cone invariant under $A^T$, we can find a strictly positive vector $e\in L$ 
such that for all $k\in\NN$ we have $(A^T)^k e \in L \subseteq \RR^n_+$. 
If $k$ is large enough, we have $A_1^{k-1}\ge \frac 12E_1$ and therefore $B_k\le A_1^{k-1}B\le \frac 12E_1 B$. 
Writing $e=e_1\oplus e_2$ with respect to the given decomposition, we get
$B_k^T e_1 \le \frac 12 (E_1 B)^T e_1 = \frac 12 B^T E_1^T e_1$. Since 
the vector $B^T E_1^T e_1$ has at least one strictly negative component, the same holds for $B_k^T e_1$.
Since $\lim_{k\to \infty} A_2^k=0$, there is some power $k$ such that the vector 
$(A^T)^k e =((A_1^T)^k e_1)\oplus (B_k^T e_1 +(A_2^T)^k e_2)$ has at least one strictly negative component. 
This is a contradiction with $(A^T)^k e \in L \subseteq \RR^n_+$.

(3) $\rho(A_1)<\rho(A_2)=1$: This case can be handled in a way similar to the case (2); we get the contradiction with the assumption that $K$ 
is a proper cone invariant under $A$. 
\qed
\end{pf}

The next example shows that in Proposition \ref{obstruction} none of the cones $K$ and $L$ can be omitted.

\begin{exm} The proper cone $K=\{(x,y)\ | \ x\ge y\ge 0\}\subset \RR_+^2$ is invariant under the matrix
$$A=\left[\matrix {1  &   -1    \cr
                   0    &  0     
}\right],$$
which is diagonally similar to a nonnegative matrix, but it is not nonnegative itself. Therefore, the cone $L$ 
cannot be omitted in Proposition \ref{obstruction}.
By duality, the cone $K$ cannot be omitted as well.
\qed
\end{exm}

The following is the main result of the paper.

\begin{thm}
\label{semi n>2} 
Let $\cS \subset M_n(\CC)$ be an irreducible semigroup such that each member of $\cS$ is diagonally similar to
a nonnegative matrix. Suppose that every member of rank at least $2$ is either indecomposable or $1$-decomposable.  
Then $\cS$ is (simultaneously) diagonally similar to a semigroup of nonnegative matrices. 
\end{thm}

\begin{pf}
By Lemma \ref{closurelemma}, we can assume that $\cS=\overline{\RR_+\cS}$. Then, by Lemma \ref{backtoreality}, we can assume that
$\cS\subset M_n(\RR)$ and that there are spanning sets $X,Y\subseteq \RR_+^{n}$ such that $\cS_1=XY^T$.
We can also assume that $X=\RR_+ X$ and $Y=\RR_+ Y$.  
Denote by $\widehat{X}$ and $\widehat{Y}$ the cones generated by $X$ and $Y$, respectively. 
Since $X$ and $Y$ are spanning sets, the cones $\widehat{X},\widehat{Y} \subseteq \RR^n_+$ are proper. 
Choose any member $A\in \cS$ of rank at least $2$. Then,  for all $x\in X$ and $y\in Y$,
the matrices $Axy^T=(Ax)y^T$ and $xy^TA=x(A^Ty)^T$ belong to $\cS_1=XY^T$. 
It follows that $Ax\in X$ and $A^Ty \in Y$, and therefore the proper cone $\widehat{X}$ is invariant under $A$, 
while the proper cone $\widehat{Y}$ is invariant under $A^T$. 
Since the matrix $A$ is either indecomposable or $1$-decomposable, we now apply either 
Lemma \ref{indecompobstruction} or Proposition \ref{obstruction} to conclude that $A$ is nonnegative. This completes the proof.
\qed 
\end{pf}

\begin{cor}
\label{semi n=2} 
Let $\cS \subset M_2(\CC)$ be an irreducible semigroup such that each member of $\cS$ is diagonally similar to
a nonnegative matrix. Then $\cS$ is (simultaneously) diagonally similar to a semigroup of nonnegative matrices. 
\end{cor}

We conclude the paper with the following example showing that the (in)decomposability assumptions in Proposition \ref{obstruction} 
and Theorem \ref{semi n>2} cannot be omitted.
 
\begin{exm}
Define the matrix
$$ A_3 = \left[\matrix {1 & 0 & 1    \cr
                      0 & 1 & -1   \cr
                      0 & 0 & 0                 
                     }\right] $$
and the proper cones $K_3=\{(x,y,z)\in\RR^3\ | \ x\ge 0\,,\ y\ge z\ge 0\}\subset\RR^3_+$ and $L_3=\{(x,y,z)\in\RR^3\ | \ x\ge y\ge 0\,,z\ge 0\}\subset\RR^3_+$.
It is easy to see that $K_3$ is invariant under $A_3$, while $L_3$ is invariant under $A_3^T$. For $n\ge 3$ we define the proper cones $K_n=K_3\oplus \RR^{n-3}_+$
and $L_n=L_3\oplus \RR^{n-3}_+$. Now we define an irreducible semigroup $\cS_1=K_n L_n^T$, consisting of matrices of rank at most $1$. We extend the matrix $A_3$ with a zero block to get a matrix $A_n=A_3\oplus 0 \in M_n(\RR)$. 
As $K_3$ is invariant under $A_3$ and $L_3$ is invariant under $A_3^T$, it is clear that 
the cones $K_n$ and $L_n$ are invariant under $A_n$ and $A_n^T$, respectively. 
Since $A_n^2 = A_n$, $\cS=\cS_1 \cup \{A_n\}$ is an  irreducible semigroup in which each member is diagonally similar to a nonnegative matrix, while the whole semigroup is not diagonally similar to a semigroup of nonnegative matrices.
\qed
\end{exm}

\end{document}